\documentclass[pdflatex,sn-mathphys-num]{sn-jnl}

\usepackage{graphicx}%
\usepackage{multirow}%
\usepackage{amsmath,amssymb,amsfonts}%
\usepackage{mathrsfs}%
\usepackage{xcolor}%
\usepackage{textcomp}%
\usepackage{manyfoot}%
\usepackage{booktabs}%
\usepackage{algorithm}%
\usepackage{algorithmicx}%
\usepackage{algpseudocode}%
\usepackage{listings}%
\usepackage{enumerate}%
\usepackage{anyfontsize}%

\theoremstyle{thmstyleone}%
\newtheorem{theorem}{Theorem}[section]
\newtheorem{lemma}[theorem]{Lemma}
\newtheorem{proposition}[theorem]{Proposition}
\newtheorem{corollary}[theorem]{Corollary}

\theoremstyle{thmstyletwo}%
\newtheorem{example}[theorem]{Example}
\newtheorem{remark}[theorem]{Remark}

\theoremstyle{thmstylethree}%

\begin{document}

\title[On Two Dimensional Flat Hessian Potentials]{On Two Dimensional Flat Hessian Potentials}

\author*{\fnm{Hanwen} \sur{Liu}}\email{hanwen.liu@warwick.ac.uk}

\affil{\orgdiv{Mathematics Institute}, \orgname{University of Warwick}, \orgaddress{\city{Coventry}, \postcode{CV4 7AL}, \country{UK}}}

\abstract{
A Riemannian metric is termed a Hessian metric if in some coordinate system it can be locally represented as the Hessian quadratic form of some locally defined smooth potential function. Under very mild extra technical conditions, we first theoretically describe the potentials of flat Hessian metrics on surfaces, and then construct these potentials explicitly using methods from integrable systems.
}


\maketitle

\section{Introduction}

Let $M$ be a differentiable manifold of dimension $n$. Recall that, for a Riemannian metric $h$ on $M$, a coordinate map $u\colon U\rightarrow\mathbb{R}^n$ of $M$ at the vicinity of a point $p\in M$ is termed a Hesse coordinate of $(M,h)$ at $p$, if and only if $$h|_U=\sum_{i=1}^n\sum_{j=1}^n\frac{\partial^2\Phi}{\partial u^i\partial u^j}d u^i\otimes d u^j$$ for some $\Phi\in C^\infty(U)$, and in this case, the smooth function $\Phi$ is termed a potential of $h$ in Hesse coordinate $u\colon U\rightarrow\mathbb{R}^n$.
A Riemannian metric $h$ on $M$ is termed a Hessian metric, if for every point $p\in M$ there exists a Hesse coordinate of $(M,h)$ at $p$.

Specifically, the following construction provides Hessian metrics: Let $f$ be a real-valued function of two real variables such that its Hessian matrix is everywhere positive definite. Then $$g:=f_{xx}dx^2+2f_{xy}dxdy+f_{yy}dy^2$$ is a Hessian metric, and $f$ is its potential function.

In this article, using methods from mathematical physics, we study thoroughly potential $f=f(x,y)$ of which Hessian quadratic form as a Riemannian metric is flat.

\section{Preliminary Results}

We shall first establish the curvature formula for Hessian surfaces:

\begin{proposition}\label{Gaußian_curvature_general}
Let $f$ be a real-valued smooth function defined on a domain in $\mathbb{R}^2$. Suppose that the Hessian matrix of $f$ is everywhere positive definite, so that $$g:=f_{xx}dx^2+2f_{xy}dxdy+f_{yy}dy^2$$ is a Riemannian metric. Then Gaußian curvature $K$ of the metric $g$ satisfies 
$$
K=-\frac{\left|\begin{array}{ccc}f_{xx} & f_{xxx} & f_{xxy} \\ f_{xy} & f_{xxy} & f_{xyy} \\ f_{yy} & f_{xyy} & f_{yyy} \end{array}\right|}{4\left|\begin{array}{cc}
f_{xx} & f_{xy} \\
f_{xy} & f_{yy}
\end{array}\right|^2}=-\frac{f_{xx}\{f_{xy},f_{yy}\}+f_{xy}\{f_{yy},f_{xx}\}+f_{yy}\{f_{xx},f_{xy}\}}{4(f_{xx}f_{yy}-f_{xy}f_{xy})^2}
$$
where $\{\cdot,\cdot\}$ is the standard Poisson bracket normalized by $\{x,y\}=1$.
\end{proposition}
\begin{proof}
Apply Brioschi's formula to Hessian metric.
\end{proof}

For flat Hessian metrics, Brioschi's formula simplifies further. 

\begin{corollary}\label{Gaußian_curvature}
Let $f$ be a real-valued smooth function defined on a domain in $\mathbb{R}^2$. Suppose that the Hessian matrix of $f$ is everywhere positive definite, so that $$g:=f_{xx}dx^2+2f_{xy}dxdy+f_{yy}dy^2$$ is a Riemannian metric. Then
$$
f_{xx}\{f_{xy},f_{yy}\}+f_{xy}\{f_{yy},f_{xx}\}+f_{yy}\{f_{xx},f_{xy}\}=\left|\begin{array}{ccc}f_{xx} & f_{xxx} & f_{xxy} \\ f_{xy} & f_{xxy} & f_{xyy} \\ f_{yy} & f_{xyy} & f_{yyy} \end{array}\right|=0
$$
if and only if $g$ is flat.
\end{corollary}
\begin{proof}
This is an immediate consequence of Proposition~\ref{Gaußian_curvature_general}.
\end{proof}

Although, as observed in \cite{4}, in dimension greater than $2$, not all Riemannian metrics are Hessian, it is proved in \cite{3} that all Riemannian metrics on surfaces are Hessian.

However, Hessian potential functions for a Riemannian metric on a surface, while exist, are in general not unique. For instance, as we shall see latter, besides the half norm square, the flat metric tensor of the Euclidean plane may well be the Hessian quadratic form of many other potential functions. 

Here, we shall provide a convenient equivalent condition for a real-valued smooth function of two real variables to be a Hessian potential of a flat metric.

\begin{lemma}\label{Homogeneous_functions}
Let $\Omega:=\{(u,v)\in\mathbb{R}^2:u^2+v^2<1\}$ be the unit disk, and $\textup{\textbf{x}}\colon \Omega\rightarrow\mathbb{R}^3$ a parametrization of a regular surface $S$ in $\mathbb{R}^3$. Suppose that $(0,0,0)\notin S$ and $\langle\textup{\textbf{x}},\textup{\textbf{x}}_u\times\textup{\textbf{x}}_v\rangle=0$. Then, there exists a homogeneous function $P\in C^\infty(\mathbb{R}^3-\{0\})$ such that $P(\textup{\textbf{x}})\equiv0$ holds at the vicinity of $(0,0)\in\Omega$.
\end{lemma}
\begin{proof}
This is a well-known fact in classical differential geometry in $\mathbb{R}^3$. We nevertheless give a sketch of proof.

Consider the smooth mapping $\textbf{F}:=\|\textbf{x}\|^{-1}\textbf{x}$. Since $\langle\textup{\textbf{x}},\textup{\textbf{x}}_u\times\textup{\textbf{x}}_v\rangle=0$, a straightforward computation shows that $\textbf{F}_u\times \textbf{F}_v\equiv\textbf{0}$. By the implicit function theorem, there exist a sufficiently small neighborhood $U$ of $(0,0)\in\Omega$ and a smooth function $$\phi\colon\{(x,y,z)\in\mathbb{R}^3:x^2+y^2+z^2=1\}\rightarrow\mathbb{R}$$ such that $\textbf{F}(U)$ is contained in the zero loci of $\phi$. Now, define homogeneous function $P$ of degree zero via $$P(x,y,z):=\phi(x/\rho,y/\rho,z/\rho)$$ where $\rho=\sqrt{x^2+y^2+z^2}$. Then, by the very construction, we have that $P(\textbf{x})\equiv0$ holds at the vicinity of $(0,0)\in\Omega$.
\end{proof}

\begin{proposition}\label{2D_flat}
Let $f$ be a real-valued smooth function defined on the unit disk $\Omega$ in $\mathbb{R}^2$. Suppose that the Hessian matrix of $f$ is everywhere positive definite, so that $$g:=f_{xx}dx^2+2f_{xy}dxdy+f_{yy}dy^2$$ is a Riemannian metric. Assume further that 
$$
\{f_{xx},f_{xy}\}^2+\{f_{xy},f_{yy}\}^2+\{f_{yy},f_{xx}\}^2>0.
$$
Then, the metric $g$ is flat if and only if there exists a homogeneous function $P\in C^\infty(\mathbb{R}^3-\{0\})$ such that $P(f_{xx},f_{xy},f_{yy})\equiv0$ holds at the vicinity of $(0,0)\in\Omega$.
\end{proposition}
\begin{proof}
Combine Corollary~\ref{Gaußian_curvature} with Lemma~\ref{Homogeneous_functions}.
\end{proof}

\section{The Main Construction}

Based on the results from section 2, to find two dimensional flat Hessian potentials, we shall search for all real-valued smooth functions in two real variables that solve the following PDE problem:
\begin{equation}\label{main_equation}
\left\{\begin{array}{l}
f_{xx}\{f_{xy},f_{yy}\}+f_{xy}\{f_{yy},f_{xx}\}+f_{yy}\{f_{xx},f_{xy}\}=0,\\f_{xx}+f_{yy}>0,\\f_{xx}f_{yy}-f_{xy}f_{xy}>0,\end{array}\right.
\end{equation}
where we also add two additional constraints
\begin{equation}\label{conditions}
\left\{\begin{array}{l}
\{f_{xx}+f_{yy},f_{xy}\}^2>0,\\ \{f_{xx},f_{xy}\}^2+\{f_{yy},f_{xy}\}^2>0,\end{array}\right.
\end{equation}
which are mild technical conditions imposed to prevent the solutions from degeneracy.

Let $f$ be a solution to equation~(\ref{main_equation}) subject to condition~(\ref{conditions}). Then, in particular, the Hessian matrix of $f$ is everywhere positive definite. For simplicity, write $(E,F,G):=(f_{xx},f_{xy},f_{yy})$. We define $u:=(E+G)^{-1}F$ and define $v:=\operatorname{log}(E+G)$. 

Since $\{E+G,F\}^2>0$, straightforward computation yields that the Jacobian $\partial(u,v)/\partial(x,y)$ is everywhere non-singular. Moreover, since $$\{E,F\}^2+\{G,F\}^2>0,$$ without loss of generality we may assume that the Jacobian $\partial(E,F)/\partial(x,y)$ is everywhere non-singular. By Proposition~\ref{2D_flat}, there exists a real-valued smooth function $\varphi$ of one real variable such that 
$$
\left\{\begin{array}{l}
E=\varphi(u)e^v,\\F=ue^v,\\G=(1-\varphi(u))e^v.
\end{array}\right.
$$
The integrability conditions $E_y=F_x$ and $F_y=G_x$ then rearrange into the following hydrodynamic system
\begin{equation}\label{hydrodynamic}
\left[\begin{array}{cc}u_y\\v_y\end{array}\right]=\frac{1}{D(u)}\left[\begin{array}{cc}u+\varphi(u)\varphi'(u) & u^2+\varphi(u)^2-\varphi(u) \\ -\varphi'(u)^2-1 & -u-\varphi(u)\varphi'(u)+\varphi'(u) \end{array}\right]\left[\begin{array}{cc}u_x\\v_x\end{array}\right]
\end{equation}
where here $D(u)=\varphi'(u)u-\varphi(u)$ is nowhere vanishing, as it is the determinant of the Jacobian $\partial(E,F)/\partial(u,v)$. Now, straightforward computation yields that the characteristic velocities 
$$
\lambda_i:=\frac{\varphi'(u)+(-1)^i\sqrt{\varphi'(u)^2-4D(u)(1+D(u))}}{2D(u)}
$$
are real and distinct for $i=1,2$. 

Take any $i\in\{1,2\}$. Recall that, up to an additive constant, the phase function $p_i$ with characteristic velocity $\lambda_i$ is the smooth function of one real variable such that $$\frac{d p_i}{d u}=\frac{1+\lambda_i\varphi'}{\lambda_i+\varphi'},$$ and the $i$-th Riemann invariant of equation~(\ref{hydrodynamic}) is $r_i(u,v):=v+p_i(u)$. The method of hodograph transformation then brings equation~(\ref{hydrodynamic}) into a system of linear equations
\begin{equation}\label{linear system}
\frac{\partial x}{\partial r_i}+\lambda_i\frac{\partial y}{\partial r_i}=0
\end{equation}
for $i=1,2$. Define the conformal coordinates $\theta:=(r_1+r_2)/2$ and $t:=(r_1-r_2)/2$. Since $2t=p_1(u)-p_2(u)$ and straightforward computation yields that
$$\frac{d p_1}{d u}>\frac{d p_2}{d u},$$ by the inverse function theorem, $u$ is a univariate function of time variable $t$. Therefore, the differential $1$-form $$\Gamma(t)dt:=\frac{2d\lambda_1}{\lambda_1-\lambda_2}$$ is well-defined. Also, up to a multiplicative constant, let $\mu$ be the positive smooth function of one real variable satisfying $2\dot{\mu}+\Gamma\mu=0$, and denote $\Psi:=\mu^2y$. Then the equation~(\ref{linear system}) reduces to the Klein-Gordon relativistic wave equation 
\begin{equation}\label{Klein-Gordon}
\ddot{\Psi}-\partial^2_\theta\Psi+V(t)\Psi=0
\end{equation}
where the potential is $V=\dot{\Gamma}-\Gamma^2$. The equation~(\ref{Klein-Gordon}) can be solved by separating variables to transform it into a Sturm-Liouville eigenvalue problem, namely, the general solution of equation~(\ref{Klein-Gordon}) is the superposition of particular solutions of the form
$$
\Psi(t,\theta)=A\operatorname{cos}(k\theta)\psi_k(t)+B\operatorname{sin}(k\theta)\psi_k(t)
$$
for some real numbers $A,B,k\in\mathbb{R}$, where $\psi_k$ is a solution of the time-independent Schrödinger equation 
$$
-\ddot{\psi}+V\psi=k^2\psi
$$
of total energy $k^2\geq0$.

Conversely, once we obtain a solution $\Psi$ of equation~(\ref{Klein-Gordon}), we immediately have that $y:=\mu^{-2}\Psi$, and we can recover $x$ via integrating the $1$-form $$dx=-\sum_{i=1}^2\lambda_i\frac{\partial y}{\partial r_i}dr_i$$
which fixes the Hesse coordinate. Then, we may invert the Hesse coordinate to obtain $u=u(x,y)$ and $v=v(x,y)$. The Hessian metric is then 
$$
H:=\left[\begin{array}{cc}f_{xx} & f_{xy} \\ f_{xy} & f_{yy} \end{array}\right]=\left[\begin{array}{cc}\varphi(u)e^v & ue^v\\ ue^v & (1-\varphi(u))e^v \end{array}\right]
$$
and its potential $f$ is locally determined by the double integration of $H$ over a star-shaped neighborhood of $(0,0)\in\Omega$.

\section{More Examples}

By exploring examples, in this section we shall show that the solutions of equation~(\ref{main_equation}) form a large family. 

\begin{example}\label{2D_homogeneous_potential}
Let $f$ be a real-valued smooth homogeneous function of degree $d$ defined on a domain in $\mathbb{R}^2$ such that its Hessian matrix $H:=\operatorname{Hess}(f)$ is everywhere positive definite. Then, by Euler's theorem on homogeneous functions, we have $$(2-d)H + x\frac{\partial H}{\partial x} + y\frac{\partial H}{\partial y} = 0,$$ that is, the matrices $H,H_x,H_y$ are everywhere $\mathbb{R}$-linearly dependent. Therefore, by Corollary~\ref{Gaußian_curvature}, we have that $$g:=f_{xx}dx^2+2f_{xy}dxdy+f_{yy}dy^2$$ is a flat Riemannian metric.
\end{example}

We also suggest the following very concrete example of a Hessian potential for the flat Euclidean plane.

\begin{example}
Let $\mathbb{H}^2:=\{(x,y)\in\mathbb{R}^2:y>0\}$ be the upper half plane and define $f\in C^\infty(\mathbb{H}^2)$ via $$f(x,y):=\frac{x^2}{2y}+\frac{1}{4}\operatorname{log}(y)y.$$ Then the coordinate change $(x(r,\theta),y(r,\theta))=(r^2\theta,r^2)$ brings the Hessian quadratic form
$$
g:=f_{xx}dx^2+2f_{xy}dxdy+f_{yy}dy^2=\frac{1}{y}dx^2-\frac{2x}{y^2}dxdy+\frac{4x^2+y^2}{4y^3}dy^2
$$
of $f$ to the flat metric tensor $ds^2=dr^2+r^2d\theta^2$ of the Euclidean plane in its polar coordinates. The flatness of $g$ can be proved also by noticing that $$f_{xx}f_{xx}-4f_{xx}f_{yy}+4f_{xy}f_{xy}\equiv0$$ and then apply Proposition~\ref{2D_flat}.
\end{example}

Before we close this section, we note the uniqueness of radially symmetric flat Hessian potential.

\begin{remark}
Let $\Omega:=\{(x,y)\in\mathbb{R}^2:x^2+y^2<1\}$ be the unit disk, and $f\in C^\infty(\Omega)$ a smooth function of which Hessian matrix is everywhere positive definite, so that $$g:=f_{xx}dx^2+2f_{xy}dxdy+f_{yy}dy^2$$ is a Riemannian metric. If $f$ is radially symmetric and $g$ is flat, then a straightforward computation yields that there exists a real number $C>0$ such that $f(x,y)=C(x^2+y^2)$.
\end{remark}

\section*{Acknowledgement}

The author thanks the reviewers for various suggestions. 
This research was completed while the author was studying at the mathematics institute of the University of Warwick. The author therefore would like to thank the hospitality of the University of Warwick.

\bibliography{bibliography}


\begin{thebibliography}{2}
\ifx \bisbn   \undefined \def \bisbn  #1{ISBN #1}\fi
\ifx \binits  \undefined \def \binits#1{#1}\fi
\ifx \bauthor  \undefined \def \bauthor#1{#1}\fi
\ifx \batitle  \undefined \def \batitle#1{#1}\fi
\ifx \bjtitle  \undefined \def \bjtitle#1{#1}\fi
\ifx \bvolume  \undefined \def \bvolume#1{\textbf{#1}}\fi
\ifx \byear  \undefined \def \byear#1{#1}\fi
\ifx \bissue  \undefined \def \bissue#1{#1}\fi
\ifx \bfpage  \undefined \def \bfpage#1{#1}\fi
\ifx \blpage  \undefined \def \blpage #1{#1}\fi
\ifx \burl  \undefined \def \burl#1{\textsf{#1}}\fi
\ifx \doiurl  \undefined \def \doiurl#1{\url{https://doi.org/#1}}\fi
\ifx \betal  \undefined \def \betal{\textit{et al.}}\fi
\ifx \binstitute  \undefined \def \binstitute#1{#1}\fi
\ifx \binstitutionaled  \undefined \def \binstitutionaled#1{#1}\fi
\ifx \bctitle  \undefined \def \bctitle#1{#1}\fi
\ifx \beditor  \undefined \def \beditor#1{#1}\fi
\ifx \bpublisher  \undefined \def \bpublisher#1{#1}\fi
\ifx \bbtitle  \undefined \def \bbtitle#1{#1}\fi
\ifx \bedition  \undefined \def \bedition#1{#1}\fi
\ifx \bseriesno  \undefined \def \bseriesno#1{#1}\fi
\ifx \blocation  \undefined \def \blocation#1{#1}\fi
\ifx \bsertitle  \undefined \def \bsertitle#1{#1}\fi
\ifx \bsnm \undefined \def \bsnm#1{#1}\fi
\ifx \bsuffix \undefined \def \bsuffix#1{#1}\fi
\ifx \bparticle \undefined \def \bparticle#1{#1}\fi
\ifx \barticle \undefined \def \barticle#1{#1}\fi
\bibcommenthead
\ifx \bconfdate \undefined \def \bconfdate #1{#1}\fi
\ifx \botherref \undefined \def \botherref #1{#1}\fi
\ifx \url \undefined \def \url#1{\textsf{#1}}\fi
\ifx \bchapter \undefined \def \bchapter#1{#1}\fi
\ifx \bbook \undefined \def \bbook#1{#1}\fi
\ifx \bcomment \undefined \def \bcomment#1{#1}\fi
\ifx \oauthor \undefined \def \oauthor#1{#1}\fi
\ifx \citeauthoryear \undefined \def \citeauthoryear#1{#1}\fi
\ifx \endbibitem  \undefined \def \endbibitem {}\fi
\ifx \bconflocation  \undefined \def \bconflocation#1{#1}\fi
\ifx \arxivurl  \undefined \def \arxivurl#1{\textsf{#1}}\fi
\csname PreBibitemsHook\endcsname

\bibitem[\protect\citeauthoryear{Amari and Armstrong}{2014}]{4}
\begin{barticle}
\bauthor{\bsnm{Amari}, \binits{S.}},
\bauthor{\bsnm{Armstrong}, \binits{J.}}:
\batitle{Curvature of {H}essian manifolds}.
\bjtitle{Differential Geometry and its Applications}
\bvolume{33},
\bfpage{1}--\blpage{12}
(\byear{2014})
\end{barticle}
\endbibitem

\bibitem[\protect\citeauthoryear{Han and Wang}{2018}]{3}
\begin{barticle}
\bauthor{\bsnm{Han}, \binits{Q.}},
\bauthor{\bsnm{Wang}, \binits{G.}}:
\batitle{{H}essian surfaces and local {L}agrangian embeddings}.
\bjtitle{Annales de l'Institut Henri Poincar\'e C, Analyse non lin\'eaire}
\bvolume{35},
\bfpage{675}--\blpage{685}
(\byear{2018})
\end{barticle}
\endbibitem

\end{thebibliography}

\end{document}